\documentclass[11pt]{article}
\usepackage{amssymb,amsfonts,amsmath,latexsym,epsf,tikz,url}

\newtheorem{theorem}{Theorem}[section]

\newtheorem{corollary}[theorem]{Corollary}
\newtheorem{lemma}[theorem]{Lemma}
\newtheorem{define}[theorem]{Definition}
\newtheorem{example}[theorem]{Example}

\newcommand{\proof}{\noindent{\bf Proof.\ }}
\newcommand{\qed}{\hfill $\square$\medskip}

\begin{document}

\title{On $\mathcal{D}$-equivalence classes of some graphs}

\author{Somayeh Jahari
\and
Saeid Alikhani $^{}$\footnote{Corresponding author}
}

\date{}

\maketitle

\begin{center}
 Department of Mathematics, Yazd University, 89195-741, Yazd, Iran\\
{\tt  alikhani@yazd.ac.ir}
\end{center}


\begin{abstract}
Let $G$ be a simple graph of order $n$.
The domination polynomial of $G$ is the polynomial
$D(G, x)=\sum_{i=1}^n d(G,i) x^i$,
where $d(G,i)$ is the number of dominating sets of $G$ of size $i$. The  $n$-barbell graph $Bar_n$  with $2n$ vertices, is formed by joining two copies of a complete graph $K_n$ by a single edge.  We prove  that for every $n\geq 2$, $Bar_n$ is not $\mathcal{D}$-unique, that is, there is another non-isomorphic graph with the same domination polynomial.   More precisely, we  show that for every $n$, the $\mathcal{D}$-equivalence class of barbell graph, $[Bar_n]$, contains many graphs, which one of them is  the complement of book graph of order $n-1$, $B_{n-1}^c$. Also we present many families of graphs in  $\mathcal{D}$-equivalence class of $K_{n_1}\cup K_{n_2}\cup \cdots\cup K_{n_k}$.
\end{abstract}


\medskip
\noindent{\bf AMS Subj. Class.:} 05C60, 05C69

\section{Introduction}

All graphs in this paper are simple of finite orders, i.e., graphs are undirected with no loops or
parallel edges and with finite number of vertices.  The {\it complement} $G^c$ of a graph $G$, is a graph with the same vertex set as $G$ and with the property that two vertices are adjacent in $G^c$ if and only if they are not adjacent in $G$.
For any vertex $v\in V(G)$, the {\it open neighborhood} of $v$ is the
set $N(v)=\{u \in V (G) | uv\in E(G)\}$ and the {\it closed neighborhood} of $v$
is the set $N[v]=N(v)\cup \{v\}$. For a set $S\subseteq V(G)$, the open
neighborhood of $S$ is $N(S)=\bigcup_{v\in S} N(v)$ and the closed neighborhood of $S$
is $N[S]=N(S)\cup S$. A set $S\subseteq V(G)$ is a {\it dominating set} if $N[S]=V$, or equivalently,
every vertex in $V(G)\backslash S$ is adjacent to at least one vertex in $S$.
The {\it domination number} $\gamma(G)$, is the minimum cardinality of a dominating set in $G$.
For a detailed treatment of domination theory, the reader is referred to~\cite{domination}.

 Let ${\cal D}(G,i)$ be the family of dominating sets of a graph $G$ with cardinality $i$ and
let $d(G,i)=|{\cal D}(G,i)|$.
The {\it domination polynomial} $D(G,x)$ of $G$ is defined as
${D(G,x)=\sum_{ i=\gamma(G)}^{|V(G)|} d(G,i) x^{i}}$ (see \cite{euro,saeid1,kotek}). This polynomial is the generating polynomial for the number of dominating sets of each cardinality.

 Calculating the domination polynomial of a graph $G$ is difficult in general, as the smallest power of a non-zero term is the domination number $\gamma (G)$ of the graph, and determining whether $\gamma (G) \leq k$ is known to be NP-complete \cite{garey}. But for certain classes of graphs, we can find a closed form expression for the domination polynomial. Two graphs $G$ and $H$ are said to be {\it dominating equivalent},
or simply ${\cal D}$-equivalent, written $G\sim H$, if
$D(G,x)=D(H,x)$. It is evident that the relation $\sim$ of being
${\cal D}$-equivalence
 is an equivalence relation on the family ${\cal G}$ of graphs, and thus ${\cal G}$ is partitioned into equivalence classes,
called the {\it ${\cal D}$-equivalence classes}. Given $G\in {\cal G}$, let
\[
[G]=\{H\in {\cal G}:H\sim G\}.
\]
We call $[G]$ the equivalence class determined by $G$.
A graph $G$ is said to be dominating unique, or simply
$\mathcal{D}$-unique, if $[G] = \{G\}$ \cite{Georj}. Determining $\mathcal{D}$-equivalence class of graphs is
one of the interesting problems on equivalence classes.

 A question of recent interest concerning this equivalence relation $[\cdot]$ asks
which graphs are determined by their domination polynomial. It is known
that cycles \cite{euro} and cubic graphs of order $10$ \cite{cubic} (particularly, the Petersen
graph) are, while if $n\equiv 0 (mod\, 3)$, the paths of order $n$ are not \cite{euro}. In \cite{complete}, authors
completely described the complete $r$-partite graphs which are $\mathcal{D}$-unique. Their results in the
bipartite case, settles in the affirmative a conjecture in \cite{ghodrat}.

\medskip

 Let $n$ be any positive integer and  $Bar_n$  be $n$-barbell graph with $2n$ vertices which is formed by joining
two copies of a complete graph $K_n$ by a single edge.
 In this paper, we consider  $n$-barbell graphs and study their domination polynomials.
 We prove  that for every $n\geq 2$, $Bar_n$ is not $\mathcal{D}$-unique.   More precisely, in Section 2, we  show that for every $n$, $[Bar_n]$ contains many graphs, which one of them is $2K_n$ and another one is  the complement of book graph of order $n-1$, $B_{n-1}^c$. In Section 3, we present many graphs in $[K_{n_1}\cup K_{n_2}\cup\cdots\cup K_{n_k}]$.

\section{ $\mathcal{D}$-equivalence classes of some graphs}

 In this section, we study the $\mathcal{D}$-equivalence classes of some graphs. First we consider the domination polynomial of the  complement
of book graph.

\begin{figure}[!ht]
\hspace{5.3cm}
\includegraphics[width=4.cm,height=2.5cm]{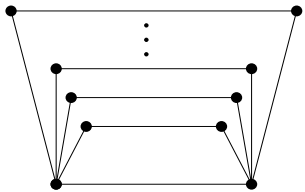}
\caption{ \label{figure6} The book graph $B_n$.}
\end{figure}

 The $n$-book graph $B_n$ can be constructed by bonding $n$ copies of the cycle graph $C_4$ along a common edge $\{u, v\}$, see Figure \ref{figure6}.

 The following theorem gives a formula for the domination polynomial of $B_n$.
\begin{theorem}{\rm \cite{jason}}
 For every $n \in \mathbb{N}$,
\[ D(B_n,x)=(x^2+2 x)^n(2x+1) + x^2(x+1)^{2n}- 2x^n.\]
\end{theorem}

  Domination polynomials, exploring the nature and location of roots of domination polynomials of book graphs has studied in \cite{jason}.  Here, we consider  the domination polynomial of the complement of the book graphs. We shall prove that the $n$-barbell graph $Bar_n$ and $B_{n-1}^c$ have the same domination polynomial.

 The Tur\'{a}n graph $T(n,r)$ is a complete multipartite graph formed by partitioning a set of $n$ vertices into $r$ subsets, with sizes as equal as possible, and connecting two vertices by an edge whenever they belong to different subsets. The graph will have $(n~ mod ~ r)$ subsets of size $\lceil\frac{n}{r}\rceil$, and $r - (n~ mod ~ r)$ subsets of size $\lfloor\frac{n}{r}\rfloor$. That is, a complete $r$-partite graph
\[ K_{\lceil\frac{n}{r}\rceil, \lceil\frac{n}{r}\rceil,\ldots,\lfloor\frac{n}{r}\rfloor,\lfloor\frac{n}{r}\rfloor}.\]

 The Tur\'{a}n graph $T(2n,n)$ can be formed by removing  a perfect matching,  $n$ edges no two of which are adjacent,  from a complete graph $K_{2n}$. As Roberts (1969) showed, this graph has boxicity exactly $n$; it is sometimes known as the Robert's graph \cite{Robert}. If $n$ couples go to a party, and each person shakes hands with every person except his or her partner, then this graph describes the set of handshakes that take place; for this reason it is also called the cocktail party graph. So, the cocktail party graph $CP(t)$ of order $2t$ is the graph with vertices $b_1, b_2, \cdots, b_{2t}$
in which each pair of distinct vertices form an edge with the exception of the pairs
$\{b_1 , b_2 \}, \{b_3 , b_4\}, \ldots, \{b_{2t- 1}, b_{2t}\}$.
 The  following result is easy to obtain.
 \begin{lemma}\label{lem}
 For every $n\in \mathbb{N}$, $D(CP(n),x) = (1+x)^{2n}-2nx-1$.
\end{lemma}

\begin{figure}[ht]
\begin{center}
\includegraphics[width=3.5in]{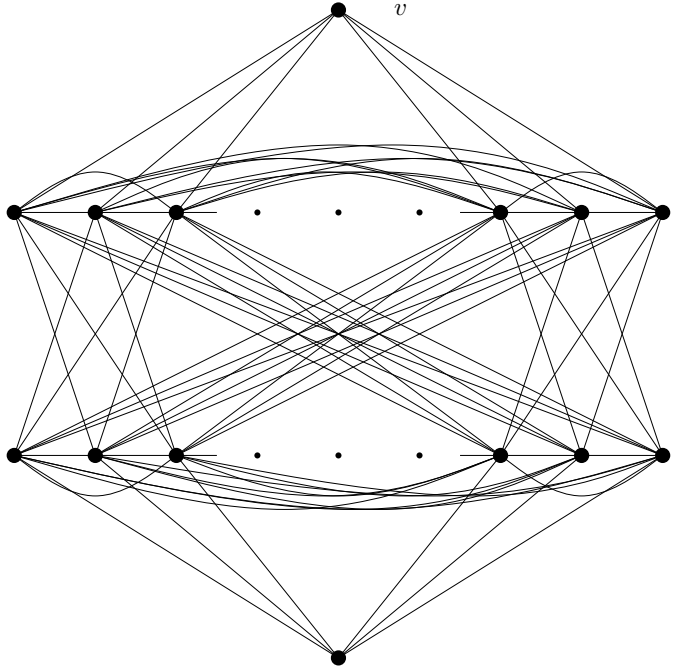}
\caption{Complement of the book graph $B_n^c$.}
\label{Cbook}
\end{center}
\end{figure}

 Figure \ref{Cbook} shows  the complement of the  book graph $B_n^c$.

\medskip
 The vertex contraction $G/u$ of a graph $G$ by a vertex $u$ is the operation under
which all vertices in $N(u)$ are joined to each other and then $u$ is deleted (see\cite{Wal}).

 The following theorem is  useful for finding the recurrence relations for the  domination polynomials  of arbitrary graphs.

\begin{theorem}\label{theorem1}{\rm \cite{SA,Kot}}
Let $G$ be a graph. For any vertex $u$ in $G$ we have
\[
D(G, x) = xD(G/u, x) + D(G - u, x) + xD(G - N[u], x) - (1 + x)p_u(G, x),
\]
where $p_u(G, x)$ is the polynomial counting the dominating sets of $G - u$ which do not contain any
vertex of $N(u)$ in $G$.
\end{theorem}

 The following theorem gives a formula for the domination polynomial of  the complement of the  book graph.

\begin{theorem}\label{theorem6}
For every $n\in \mathbb{N}$,
$$D(B_n^c,x) = ((1+x)^{n+1}-1)^2.$$
\end{theorem}

\proof
 Consider graph  $B_{n}^c$ and  vertex $v$ in the  Figure \ref{Cbook}. By Theorem \ref{theorem1}, we have:
\begin{eqnarray*}
D(B_n^c, x)&=& x D(B_n^c/v, x) + D(B_n^c - v, x) + x D(B_n^c - N[v], x) - (1 + x)p_v(B_n^c, x)\\
&=& (x +1) D(B_n^c-v, x) + xD(K_{n+1},x) - (1 + x)(D(K_{n+1},x)-(n+1)x-nx^2)\\
&=&(x +1) D(B_n^c-v, x) - D(K_{n+1},x)+x(1+x)(1+n(1+x))\\
&=&(x +1) D(B_n^c-v, x) - ((1+x)^{n+1}-1)+x(1+x)(1+n(1+x)),
\end{eqnarray*}
 where $(B_n^c/v) \simeq B_n^c-v$.

 Now,  we use Theorem \ref{theorem1} to obtain the domination polynomial of the graph $B_n^c-v$. We have
\begin{eqnarray*}
D(B_n^c-v, x)&=& x D(B_n^c-v/u, x) + D((B_n^c - v)-u, x) \\
&&+ x D((B_n^c-v) - N[u], x) - (1 + x)p_u(B_n^c-v, x).
\end{eqnarray*}
Since $(B_n^c-v/u) \simeq (B_n^c-v)-u\simeq CP(n)$ and using Lemma \ref{lem}, we have
\begin{eqnarray*}
D(B_n^c-v, x)&=& (x +1) D(CP(n), x) + x(D(K_{n},x)) - (1 + x)(D(K_{n},x)-nx)\\
&=&(x +1) D(CP(n), x) - D(K_{n},x)+nx(1+x)\\
&=&(x +1)((1+ x)^{2n} -(1 + 2nx)) - ((1+x)^n-1)+nx(1+x)\\
&=&(1+x)^n((1+x)^{n+1}-1)-nx(1+x)-x.
\end{eqnarray*}
  Consequently,
\begin{eqnarray*}
D(B_n^c, x)&=&(x +1) ((1+x)^n((1+x)^{n+1}-1)-nx(1+x)-x)\\
&& - ((1+x)^{n+1}-1)+x(1+x)(1+n(1+x))\\
&=&((1+x)^{n+1}-1)^2.
\end{eqnarray*}
\qed

\begin{figure}[!ht]
\hspace{4cm}
\includegraphics[width=5.5cm,height=2.5cm]{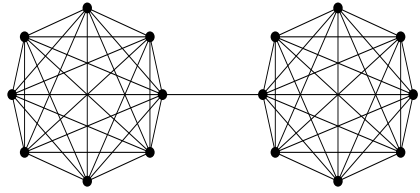}
\caption{ \label{figure1} The barbell graph of order 16, $Bar_8$.}
\end{figure}

 The $n$-barbell graph is the graph on $2n$ vertices which is formed by joining
two copies of a complete graph $K_n$ by a single edge, known as a bridge, shown in Figure \ref{figure1}. We denote this graph by $Bar_n$. For this graph, we shall calculate this domination polynomial. We need the following definition and theorems.

 An irrelevant edge is an edge $e\in E(G)$, such that $D(G, x) = D(G-e, x)$, and a vertex $v \in V(G)$ is domination-covered,
if every dominating set of $G-v$ includes at least one vertex adjacent to $v$ in $G$ \cite{Kot}.
 We need the following theorems to obtain the domination polynomial of barbell graph $Bar_n$.

\begin{theorem}{\rm\cite{Kot}}\label{dcov}
Let $G = (V,E)$ be a graph. An vertex  $ v\in V$ of $G$ is domination-covered if and only if there is a vertex  $u\in N[v]$ such that $N[u] \subseteq N[v]$.
\end{theorem}

\begin{theorem}{\rm\cite{Kot}}\label{irre}
Let $G = (V,E)$ be a graph. An edge $e = \{u, v\} \in E$ is an irrelevant edge
in $G$, if and only if $u$ and $v$ are domination-covered in $G-e$.
\end{theorem}

\begin{theorem}\label{theorem5}
For every $n\geq 2$ and $n\in \mathbb{N}$,
$$D(Bar_n,x) = ((1+x)^{n}-1)^2.$$
\end{theorem}

\proof
 Let $e$ be an  edge joining two $K_n$ in barbell graph.
By Theorem \ref{dcov} two end vertices of edge $e$ are  domination-covered in $Bar_n-e$. So, by Theorem \ref{irre}  the edge $e$ is an  irrelevant edge of $Bar_n$. Therefore
 $$D(Bar_n,x)=D(Bar_n-e,x)=D(K_n\cup K_n,x)=((1+x)^{n}-1)^2.$$
 \qed

 The following corollary is an immediate consequence of Theorems \ref{theorem6} and \ref{theorem5}.
\begin{corollary}\label{korolari1}
For each natural number $n$,  $Bar_n$ and $B_{n-1}^c$ have the same domination polynomial. More precisely,  for every $n$, $[Bar_n]\supseteq \{Bar_n,B_{n-1}^c,K_n\cup K_n\}$, and $[B_{n-1}^c]\supseteq \{Bar_n,B_{n-1}^c,K_n\cup K_n\}$.
\end{corollary}

\medskip

\begin{figure}[h]
\hspace{3cm}
\includegraphics[width=8cm,height=2.4cm]{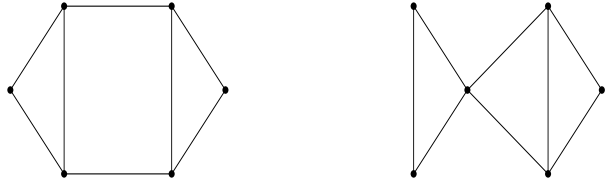}
\caption{\label{gen} Two generalized barbell graphs $Bar_{3,2}$. }
\end{figure}

 Here, we present some other families of graphs whose are in the $[Bar_n]$. Let to define the generalized barbell graphs. As we know, the $Bar_n$ is  formed by joining two copies of a complete graph $K_n$ by a single edge. We like to join two copies with more edges as follows:

 \begin{define}
 Suppose that $\{u_1,...,u_n\}$ and $\{v_1,...,v_n\}$ are the vertices of two copies of complete graph of order $n$, $K_n$ and $\mathcal{K}_n$.
  The generalized barbell graph is denoted by $Bar_{n,t}$ and is a graph with $V(Bar_{n,t})=\{u_1,...,u_n\}\cup \{v_1,...,v_n\}$ and
   $$E(Bar_{n,t})=E(K_n)\cup E(\mathcal{K}_n)\cup\big\{u_iv_j|1\leq i\leq n-1, 1\leq j\leq n-1\big\},$$
     where $\Big|\big\{u_iv_j|1\leq i\leq n-1, 1\leq j\leq n-1\big\}\Big|=t$.
\end{define}
  As examples see two non-isomorphic graphs $Bar_{3,2}$ in Figure \ref{gen}.   Notice that $B_{n-1}^c$ is one of the specific case of $B_{n,(n-1)(n-2)}$.  The left graph in Figure \ref{gen}, is $B_2^c$.

 We have the following theorem.

\begin{theorem}\label{theorem7}
For every $n\geq 3$ and $n\in \mathbb{N}$,
$$D(Bar_{n,t},x) = ((1+x)^{n}-1)^2.$$
\end{theorem}

 \proof
  We prove  this Theorem by induction on $t$. Suppose that $t = 1$, Then by Theorem \ref{theorem5}, the result holds. Assume that the result holds for $t= (n-1)^2-1$. Let $t=(n-1)^2$ and $e$ be the additional edge of $Bar_{n,t}$ to the $Bar_{n,t-1}$. By Theorem \ref{dcov} two end vertices of edge $e$ are  domination-covered in $Bar_n-e$. So, by Theorem \ref{irre}  the edge $e$ is an  irrelevant edge of $Bar_{n,t-1}$. Therefore by the induction hypothesis we have the result.
  \qed

 The following corollary is an immediate consequence of Theorems \ref{theorem5} and \ref{theorem7}.
\begin{corollary}
 For each natural number $n$ and $t\leq (n-1)^2$,  $Bar_n$ and $Bar_{n,t}$ have the same domination polynomial.
\end{corollary}

 The following example shows that, except for the generalized barbell graphs, there are other graphs in $\mathcal{D}$-equivalence classes of $Bar_n$.

\begin{figure}[h]
\hspace{4.6cm}
\includegraphics[width=6cm,height=2.4cm]{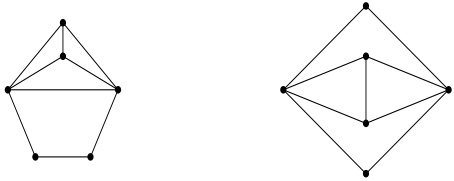}
\caption{\label{othe} Two graphs in $[Bar_3]$. }
\end{figure}

\begin{example}\label{last}
 All connected graphs in $[Bar_3]$ are the graphs $Bar_3,~Bar_{3,2},~Bar_{3,3},~Bar_{3,4}$ and two  graphs in Figure \ref{othe}.
 \end{example}

\section{Some graphs in $[K_{n_1}\cup K_{n_2}\cup\cdots\cup K_{n_k}]$}

 We observed that, for each natural number $n$ and $t\leq (n-1)^2$,  the domination polynomials of  $Bar_n$ and $Bar_{n,t}$ is $((1+x)^n-1)^2$.  In this section, we present graphs whose domination polynomials are $\prod_{i=1}^k((1+x)^{n_i}-1)$. For this purpose, we construct   families of  graphs from a path $P_k$ which we denote them by $S(G_1,G_2,...,G_k)$ in the following definition.

 \begin{define} \label{def1}
 The graph $S(G_1,G_2,...,G_k)$   is a graph which obtain from a path $P_k$ with the vertices $\{v_1, v_2, \ldots, v_k\}$, by
substituting a graph  $G_i$ of order $n_i\geq 3$,  for every vertex $v_i$ of $P_k$, such that
\begin{itemize}
\item for $i=1,k$, the graphs $G_i$ have at least one vertex of degree $n_i-1$ and other $G_i$'s have  at least two vertices of degree $n_i-1$, and
\item  in the graph $S(G_1,G_2,...,G_k)$, the end vertices of each edge $e_i$ in the path graph, $P_k$ are one vertex of degree $n_i-1$ in graphs $G_{i-1}$ and $G_i$.
\end{itemize}
\end{define}
 We have the following result for  graph $S(G_1,G_2,...,G_k)$.
\begin{theorem}
For every natural number $k\geq 2$,
$$D(S(G_1,G_2,...,G_k),x) = D(G_1,x)D(G_2,x)\ldots D(G_k,x).$$
In particular if $G_i=K_{n_i}$ and $n_i\geq 3$, then
$$D(S(K_{n_1},K_{n_2},...,K_{n_k}),x) = \prod_{i=1}^k D(K_{n_i},x)=\prod_{i=1}^k((1+x)^{n_i}-1).$$
\end{theorem}

\proof
 Let $e_i$ $(1\leq i\leq k)$ be the  edge joining $G_{i-1}$ and  $G_i$ in $S(G_1,G_2,...,G_k)$.
By Theorem \ref{dcov} two end vertices of edge $e_i$ are  domination-covered in $S(G_1,G_2,...,G_k)-e_i$. So, by Theorem \ref{irre}  every edge $e_i$ is an  irrelevant edge of $S(G_1,G_2,...,G_k)$. Therefore we have the result.
\qed

 We shall  generalize the graphs $S(G_1,G_2,...,G_k)$ in Definition \ref{def1} such that this generalized graphs and $S(G_1,G_2,...,G_k)$ have the same
 domination polynomial. Suppose that $GS_t(K_{n_1},K_{n_2},...,K_{n_k})$
be a family of graphs in the form of $S(K_{n_1},K_{n_2},...,K_{n_k})$  such that the complete graphs $K_{n_i}$ with  $V(K_{n_i})=\{u_1,...,u_{n_i}\}$ and $K_{n_{i+1}}$ with $V(K_{n_{i+1}})= \{v_1,...,v_{n_{i+1}}\}$ are joined with $t$ following  edges
 \[\big\{u_iv_j|1\leq i\leq n_{i}-1, 1\leq j\leq n_{i+1}-1\big\}.\]

Similar to the proof of the Theorem \ref{theorem7}, we have the following theorem:

\begin{theorem}
For each natural number $t$, all graphs in the family of  $GS_t(K_{n_1},K_{n_2},...,K_{n_k})$ have the same domination polynomial. More precisely,   the domination polynomial of each $H$ in $GS_t(K_{n_1},K_{n_2},...,K_{n_k})$ is equal to $\prod_{i=1}^k((1+x)^{n_i}-1)$.
\end{theorem}

\medskip

\noindent{\bf Conclusion.}
In this paper, we studied the $\mathcal{D}$-equivalence classes of barbell graphs $Bar_n$. We showed that, for each natural number $n$, $2K_n$,  $Bar_n$, $Bar_{n,t}$ and the complement of the book graph of order $n-1$, $B_{n-1}^c$ have the same domination polynomial, i.e., $[Bar_n]=[Bar_{n,t}]=[B_{n-1}^c]=[K_n\cup K_n]$.
Example \ref{last}, implies that except for these kind of graphs, there are another graphs in this class. Therefore, exact characterization of graphs in $[Bar_n]$ remains as an open problem. Also we presented many families of graphs whose are in  $[K_{n_1}\cup K_{n_2}\cup...\cup K_{n_k}]$, but similar to Example \ref{last}, there are another graphs in this class. So,  exact characterization of graphs in $[K_{n_1}\cup K_{n_2}\cup...\cup K_{n_k}]$ remains as another open problem.

\end{document}